\documentclass{amsart}

\usepackage{amssymb}

\usepackage[all]{xy}

\newtheorem{thm}{Theorem}

\newtheorem{lem}[thm]{Lemma}

\newtheorem{prop}[thm]{Proposition}

\theoremstyle{definition}
\newtheorem{defn}[thm]{Definition}

\newtheorem{say}[thm]{}

\newtheorem{rem}[thm]{Remark}          

\newtheorem*{ack}{Acknowledgments}      

\newtheorem{defn-thm}[thm]{Definition--Theorem}  
\newtheorem{defn-lem}[thm]{Definition--Lemma}  

\theoremstyle{remark}

\renewcommand{\c}[0]{{\mathbb C}}  

\renewcommand{\o}[0]{{\mathcal O}} 
\newcommand{\z}[0]{{\mathbb Z}}

\renewcommand{\a}[0]{{\mathbb A}}

\newcommand{\p}[0]{{\mathbb P}}

\newcommand{\q}[0]{{\mathbb Q}}

\newcommand{\qtq}[1]{\quad\mbox{#1}\quad}

\newcommand{\supp}[0]{\operatorname{Supp}}

\newcommand{\coker}[0]{\operatorname{coker}}

\newcommand{\sing}[0]{\operatorname{Sing}}

\newcommand{\onto}[0]{\twoheadrightarrow}

\newcommand{\diag}[0]{\operatorname{diag}}

\newcommand{\simq}[0]{\sim_{\q}}

\newcommand{\tsum}[0]{\textstyle{\sum}}

\def\into{\DOTSB\lhook\joinrel\to}

\def\loccoh#1.#2.#3.#4.{H^{#1}_{#2}(#3,#4)}

\DeclareMathAlphabet{\mathchanc}{OT1}{pzc}%
                                {m}{it}

\def\CC{\mathcal C}
\def\DD{\mathcal D}
\newcommand{\DMR}{{\mathcal{DMR}}}
\newcommand{\DR}{{\mathcal{DR}}}

\begin{document}
\bibliographystyle{amsalpha}

\title[Simple normal crossing varieties]
{Simple normal crossing varieties with
\\ prescribed dual complex}
\author{J\'anos Koll\'ar}

\maketitle

Let $W$ be a  reducible variety
with simple normal crossing singularities only.
 The combinatorial structure of  $W$ is described by
the {\it dual complex}  $\DD(W)$ whose vertices correspond to
the irreducible components $W_i\subset W$ and the positive dimensional 
cells correspond to the  {\it strata}
of $W$, that is, to  irreducible components of  intersections
$\cap_{i\in J} W_i$ for some $J\subset I$;
 see Definition \ref{snc-nc.sing.defn}.

In the papers
\cite{k-t-lc-exmp, k-fg,  k-link} 
many interesting singularities were obtained
by constructing a simple normal crossing variety $W$ whose
dual complex is topologically complicated and then
realizing $W$ as the exceptional divisor of a (partial)
resolution of a singularity.

The aim of this note is to prove the following existence result
for simplicial complexes.

\begin{thm}\label{main.main.thm} 
 Let $\CC$ be a finite simplicial complex
of dimension $\leq n$. Then there is a
smooth projective variety $Y$ of dimension $n+1$ and a
simple normal crossing divisor $D\subset Y$ 
such that  $\DD(D)\cong \CC$.

Moreover, we can achieve that either one (but not both)
of the following conditions are also satisfied.
\begin{enumerate}
\item $Y$ and all the strata of $D$ are rational or
\item $K_Y+D$ is ample.
\end{enumerate}
\end{thm}

This is stronger than the results in \cite{k-fg} in three aspects.
We realize $\CC$ up-to isomorphism
(not just up-to  homotopy equivalence) and
the dimension estimate  $\dim D\geq \dim \CC$   is optimal
since $\dim W\geq \dim \DD(W)$  for every simple normal crossing variety $W$.
Furthermore,   by construction $D$ is naturally  a 
hypersurface in a smooth projective variety, while
in earlier versions this was achieved only after
a detailed study of a resolution process in \cite[Sec.6]{k-link}.
The last property is very helpful in 
constructing singularities. It  simplifies the
method of \cite{k-t-lc-exmp} and also leads to isolated singularities.

\begin{thm}\label{main.sing.thm} 
 Let $\CC$ be a finite simplicial complex
of dimension $\leq n$. Then there is an isolated singularity
$(x\in X)$ of dimension $n+1$ and a log resolution
$\pi:Y\to X$ with exceptional divisor $E=\supp \pi^{-1}(x)$ 
such that  $\DD(E)\cong \CC$.
\end{thm}

Both of the alternatives in Theorem  \ref{main.main.thm}
are useful. The first case  (\ref{main.main.thm}.1)
can be utilized to  complete the  characterization of
dual complexes of resolutions of 
isolated rational singularities started in 
\cite{k-fg, k-link}.
 
\begin{thm}\label{main.rtl.thm} 
For a finite  simplicial complex $\CC$
the following are equivalent.
\begin{enumerate}
\item $\CC$ is connected, $\dim \CC\leq n$   and
$H_i\bigl(\CC, \q\bigr)=0$ for $i>0$.
\item There is an isolated rational singularity
$(x\in X)$ of dimension $n+1$ and a log resolution
$\pi:Y\to X$ with exceptional divisor $E=\supp \pi^{-1}(x)$ 
such that  $\DD(E)\cong \CC$.
\end{enumerate}
\end{thm}

Let $(x\in X)$ be a singularity and
$g:Y\to X$ a  resolution such that $E:=\supp g^{-1}(x)\subset Y$
is a simple normal crossing divisor.
The dual complex $\DD(E)$ depends on the choice of $Y$ but
its (simple) homotopy type does not; see Definition~\ref{dres.defn}.
This (simple) homotopy type  is denoted by
$\DR(x\in X)$. 

If $K_X$ is Cartier (or at least $\q$-Cartier)
then \cite{dkx} identifies a distinguished PL-homeomorphism class
$\DMR(x\in X)$ within the homotopy class $\DR(x\in X)$;
see Definition \ref{dres.defn.II}.
Using the alternative (\ref{main.main.thm}.2),
we see that the PL-homeomorphism classes  $\DMR(x\in X)$
can be arbitrary.

\begin{thm}\label{main.DMR.thm} 
For every finite,  connected  simplicial complex $\CC$
of dimension $\leq n$ 
there is a normal,  isolated singularity
$(x\in X)$ of dimension $n+1$ 
such that  $K_X$ is Cartier and 
$\DMR(x\in X)$ is PL-homeomorphic to $\CC$.
\end{thm}

\begin{say}[Open problems] \label{open.pr.say}
There are four main open problems in connection with
dual complexes of varieties and of singularities.
\medskip

(\ref{open.pr.say}.1) It is interesting to understand the possible
relation between the dual complex $\CC(X)$ and the canonical
class $K_X$.

If $K_X$ is ample then $\CC(X)$ can be arbitrary
by (\ref{main.main.thm}.2).

 If $-K_X$ is ample then $\CC(X)$ is a simplex by 
\cite{2012arXiv1206.1994F, 2012arXiv1206.2475F}.

The really interesting case seems to be when $K_X\sim 0$
or, more generally, when $K_X\simq 0$.
These lead to log canonical singularities and their
structure is not yet understood.

By \cite{k-source}  (see also \cite[Sec.4.4]{kk-singbook})
in this case either $\dim \CC(X)=1$ or 
 $\CC(X)$ is a normal pseudomanifold. In particular, every 
codimension 1 cell is the face of exactly 2 
maximal dimensional cells.

In all examples that I know of, $\CC(X)$ is an orbifold.
\medskip

(\ref{open.pr.say}.2) In many cases, the most economical realization
of a PL-homeomor\-phism type is given not by a 
simplicial complex but by a $\Delta$-complex  in the terminology of
 \cite[p.534]{hatcher}. 
Thus it would be quite useful to know that every
$\Delta$-complex can be realized as the dual complex of
 a normal crossing variety. This would give the optimal answer
along the lines proposed by \cite{simp}.
See \cite{2012arXiv1201.3129K}  for very interesting related results.
\medskip

(\ref{open.pr.say}.3) Usually one can not satisfy
both of the alternatives in Theorem  \ref{main.main.thm}
simultaneously. Indeed, if $v_{n-1}\in \DD(D)$ is an $(n-1)$-cell
that is the face of 1 or 2 $n$-cells then the corresponding
stratum is a smooth curve $C\subset D$ with 1 or 2 smaller strata
on it. Thus if $C$ is rational then $(K_D\cdot C)=-1$ or $(K_D\cdot C)=0$.
This is, however, the only obstruction that I know.

It is possible that if $\CC$ is a pseudomanifold
then there is a simple normal crossing variety $X$ 
such that  all the strata of $X$ are rational,
$K_X$ is nef and 
$\DD(X)\cong \CC$.

\medskip

(\ref{open.pr.say}.4)
A shortcoming of our construction is that it is not local.
Let $X=\cup_{i\in I} X_i$ be an snc variety. The link of a vertex
$v_i\in \DD(X)$  is determined by those $X_j$ 
for which $X_j\cap X_i\neq \emptyset$. Conversely, 
give a complex $\CC$ one could hope to that there is a
procedure to construct an  snc variety $X=\cup_{i\in I} X_i$
such that $\DD(D)\cong \CC$ and each
 $X_i$ is determined by the link of $v_i\in \CC$. 
Such a procedure exists if $\dim \CC=2$; see \cite{k-t-lc-exmp}.
The study of parasitic intersections in \cite{k-fg}
suggests that this may be too much to hope, but leaves 
open the possibility that $X_i$ should depend
only on the $(\dim \CC-1)$-times iterated link of $v_i$

One clearly needs a more local construction 
if one hopes to control the canonical class $K_X$.
\end{say}

\section{Dual complexes}

Since our main interest is in the topological aspects of these questions,
we work with varieties over $\c$. However, everything applies
to varieties over arbitrary fields. See \cite{kk-singbook}
for the more delicate aspects that appear when the
base field is not algebraically closed and has positive characteristic.

\begin{defn}[Simple normal crossing varieties]

Let $W$ be a variety
with irreducible components $\{W_i: i\in I\}$.
We say that $W$ is a {\it simple normal crossing} variety
(abbreviated as {\it snc}) if the $W_i$ are smooth and 
 every point $p\in W$ has an open (Euclidean) neighborhood $p\in U_p\subset W$
and an  embedding $U_p\into \c^{n+1}_p$ 
such that  the image of $U_p$ is an open subset of
the union of coordinate hyperplanes $ (z_1\cdots z_{n+1}=0)$.

A {\it stratum}
of $W$ is any irreducible component of an intersection
$\cap_{i\in J} W_i$ for some $J\subset I$.

The intersection of any two strata is also a union of 
lower dimensional strata, hence the union of all minimal strata
is a smooth (reducible) subvariety of $W$.
\end{defn}

\begin{defn}[Dual complex]\label{snc-nc.sing.defn}
The combinatorics of a simple normal crossing variety
$W$ is encoded by a cell complex $\DD(W)$
 whose vertices are labeled by the
 irreducible components of $W$ and for every stratum
$Z\subset \cap_{i\in J} W_i$ we attach a $(|J|-1)$-dimensional cell.
Note that for any $j\in J$ there is a unique  irreducible component
of $\cap_{i\in J\setminus\{j\}} W_i$ that contains $Z$;
this specifies the attaching map.
$\DD(W)$ is called the
{\it dual complex} or  {\it  dual graph} of $W$. 
(Although $\DD(W)$ is not a simplicial complex in general, it is a
{\it regular cell complex}; cf.\
\cite[p.534]{hatcher}.)

The minimal strata of $W$ correspond to the  maximal cells of $\DD(W)$.

It is clear that $\dim \DD(W)\leq \dim W$ and $\DD(W) $ is connected iff
$W$ is. More generally, the cohomology groups $H^i\bigl(\DD(W), \q\bigr)$
can be identified with the weight 0 part of the
Hodge structure on  $H^i\bigl(W, \q\bigr)$; see  \cite{2011arXiv1102.4370A}.
One can think of $\DD(W)$ as the combinatorial part of the
topology of $W$. 
\end{defn}

\begin{defn}[Simple normal crossing pairs]
Let $X$ be a variety over $\c$ and $D\subset X$ a divisor.
We say that $(X, D)$ is a  {\it simple normal crossing} pair
if 
\begin{enumerate}
\item $X$ is smooth, 
\item the irreducible components $D_i\subset D$ are smooth and
\item  every point $x\in X$ has an  Euclidean open neighborhood 
$x\in U_x\subset X$
with local coordinates $z_1,\dots, z_n$ such that
$D\cap U_x= (z_1\cdots z_r=0)$ for some $r$.
\end{enumerate}

If $X$ is smooth then  $(X, D)$ is a   simple normal crossing pair iff 
$D$ is a  simple normal crossing variety.
(In Lemma \ref{elem.dlt.lem} the key point  is to understand what happens when
$D$ is a  simple normal crossing variety but $X$ is not smooth.)

One can think of a  simple normal crossing variety
$D=\cup_{i\in I}D_i$ as being glued together
from the simple normal crossing pairs
$$
\bigl(D_i, \tsum_{j\neq i} D_j|_{D_i}\bigr).
$$

We say that  $(X, D)$ is a  {\it  normal crossing} pair
if it satisfies conditions (1) and (3). 
The difference between the two variants is technical but
important for several of our results.
\end{defn}

\begin{say}[Blowing up strata]\label{bu.strata.defn}
Let $W$ be a simple normal crossing variety and $Z\subset W$ a stratum.
Then the blow-up  $B_ZW$ is again a simple normal crossing variety.
To see this, take a local chart where
$$
(Z\subset W)=\bigl((x_1=\dots=x_r=0)\subset (x_1\cdots x_n=0)\bigr)
\subset \c^{n+m}.
$$
The blow-up  $B_Z\c^{n+m}$ is defined by
$$
 B_Z\c^{n+m}=\bigl(x_is_j=x_js_i: 1\leq i, j\leq r\bigr)
\subset 
\c^{n+m}_{\mathbf x}\times \p^{r-1}_{\mathbf s}
$$
and the blow-up  $B_ZW$ is defined by
$$
B_ZW=(s_1\cdots s_r \cdot x_{r+1}\cdots x_n=0)\subset B_Z\c^{n+m}.
$$
In  a typical local chart $U_r:=(s_r\neq 0)$ we can use new coordinates
$$
\bigl(x'_1=s_1/s_r, \dots, x'_{r-1}=s_{r-1}/s_r, s_r, 
x_{r+1}, \dots, x_{n+m}\bigr)
$$
to get that
$$
B_ZW\cap U_r=(x'_1\cdots x'_{r-1}\cdot x_{r+1}\cdots x_n=0).
$$
Let $v_Z\subset \DD(W)$ be the cell corresponding to $Z$.
We see from the above computations that
$\DD\bigl(B_ZW\bigr)$ is obtained from $\DD(W)$
by removing all the cells that have $v_Z$ as a face.
(That is, removing the star of $v_Z$.)
In particular, if $Z$ is a  minimal stratum then
we remove the cell $v_Z$.

Although we do not need it, it is useful to know that
there is another way to blow up a stratum.
If $(X, D)$ is a simple normal crossing pair and
$Z\subset D$ a stratum, one can consider
the blow-up $\pi:B_ZX\to X$ with exceptional divisor $E_Z$. Then
$B_ZD\cong \pi^{-1}_*D$ and 
$\bigl(B_ZX, E_Z+B_ZD\bigr)$
is again a simple normal crossing pair. It is not hard to see that
$\DD\bigl(E_Z+B_ZD\bigr)$ is obtained as the
stellar subdivison of $\DD(D)$ with vertex $v_Z$.

\end{say}

\begin{defn}[Dual complex of a singularity I]\label{dres.defn}
 Let $X$ be a  variety and
$x\in X$ a point. Choose a resolution of
 singularities  $\pi:Y\to X$ such that
$E_x:=\pi^{-1}(x)\subset Y$ is a simple normal crossing divisor.
Thus it has a dual  complex  $\DD(E_x)$.

The dual graph of a normal surface singularity has a long history
but systematic investigations in higher dimensions 
were started only recently; see
\cite{MR2320738, MR2399025, payne11, 2011arXiv1102.4370A}. 
It is proved in these papers
 that  the homotopy type  (even the simple homotopy type) of 
$\DD(E_x)$
is  independent of the resolution $Y\to X$. We denote it by
$\DR(x\in X)$.  
($\DR$ stands for the Dual complex of a Resolution.)

It is clear that $\dim \DR(x\in X)\leq \dim X-1$ and
$\DR(x\in X)$ is connected if $X$ is normal.
\end{defn}

\begin{defn}[Dual complex of a singularity II]\label{dres.defn.II}
It was observed in \cite{dkx} that the
dual complex of a divisor $D\subset Y$ is defined whenever
$(Y,D)$ is {\it divisorial log terminal} (abbreviated as {\it dlt}); see
\cite[2.37]{km-book} for the definition of dlt. 
If $(x\in X)$ is an isolated singularity and
$K_X$ is $\q$-Cartier then  there is a proper birational morphism
  $\pi:Y\to X$ such that
$(Y,E_Y)$ is dlt and $K_Y+E_Y$ is $\pi$-nef where
$E_Y=\supp \pi^{-1}(x)$.  
(Such a morphism is called a dlt modification.)
It is proved in 
 \cite{dkx} that 
\begin{enumerate}
\item $\DD(E_Y)$ 
is simple homotopy equivalent to $\DR(x\in X) $ and
\item up-to PL-homeomor\-phism, $\DD(E_Y)$ 
does not depend on the choice of
 $Y$. 
\end{enumerate}
Thus, for such singularities,   $\DD(E_Y)$  defines a distinguished
 PL-homeomorphism class, denoted by
$\DMR(x\in X)$,  within the simple homotopy class $\DR(x\in X)$.
($\DMR$ stands for the Dual complex of a Minimal partial Resolution.)

The hard part of  \cite{dkx} is the proof of (\ref{dres.defn.II}.1),
but in our examples we check this directly in Lemma \ref{elem.dlt.lem}.
\end{defn}

\section{Construction of simple normal crossing varieties}

In \cite{k-fg} we constructed snc varieties
$W=\cup_{i\in I} W_i$ by first obtaining the irreducible components
$W_i$ and the gluing them together to get $W$. 
This had two technical disadvantages.
First, the projectivity of the resulting $W$ was hard to control
and there was no natural embedding of $W$ into any smooth variety.


The  process  turns out to be much easier 
if we break the symmetry on the algebraic side 
but preserve it topologically.
Given $\CC$, first we construct an snc variety $X_0$ such 
that  $\DD(X_0)$ has the same dimension and the same vertex set as $\CC$
but the maximum possible number of positive dimensional cells.
 Then we remove the excess positive dimensional  cells one dimension at a time.

A disadvantage is that we seem to lose information
about finite group actions, while 
 the construction of \cite{k-fg}  preserved
symmetries.

The key step is the following.

\begin{lem}\label{snc.bu.induct.lem}  Let $\CC$ be a finite simplicial complex
of dimension $n$ with  vertices $\{v_i:i\in I\}$.
Fix $0\leq r<n$.
Let $X_r$ be  an snc variety 
such that $ \DD(X_r)$ is a simplicial complex and 
$\tau_r:\CC\into \DD(X_r)$ an  embedding
that is  a bijection on the vertices
and on cells of dimension $\geq n-r+1$.

Let  $J_{n-r+1}\subset \binom{I}{n-r+1}$ be those  $(n-r+1)$-element subsets
that do not span an $(n-r)$-cell in $\CC$ and
$Z_r\subset X_r$ the union of the strata 
that correspond to elements of $\tau_r\bigl(J_{n-r+1}\bigr)$.  Then
\begin{enumerate}
\item $Z_r$ is  smooth,
\item $X_{r+1}:=B_{Z_r}X_r $, the blow-up of $Z_r\subset X_r$,  is an snc variety,
\item $\DD(X_{r+1})\subset \DD(X_r)$ and
\item the restriction of $\tau_r$ gives an embedding
$\tau_{r+1}:\CC\into \DD(X_{r+1})$
which is  a bijection on the vertices
and on cells of dimension $\geq n-r$.
\end{enumerate}
\end{lem}

Proof. 
If $s_{n-r+1}\subset I$ does not 
span an $(n-r)$-cell in $\CC$ then no subset of $I$
containing $s_{n-r+1}$ spans a cell in $\CC$. 
Thus either $s_{n-r+1}$ does not span a cell in $\DD(X_r)$
or it spans a maximal cell in $\DD(X_r)$.
Since maximal cells correspond to minimal strata,
$Z_r$ is a union of minimal strata.
The minimal strata are disjoint from each other, hence
$Z_r$ is  smooth and of pure dimension $r$.

Thus, by Paragraph \ref{bu.strata.defn},  $X_{r+1}=B_{Z_r}X_r $ is an snc variety
and $\DD(X_{r+1})$ is obtained from $\DD(X_r)$
by removing all $(n-r)$-cells that do not 
correspond to a cell in $\CC$. \qed
\medskip

The following yields  a proof of
(\ref{main.main.thm}.1).

\begin{prop}\label{pf.main.main.1}
 Let $\CC$ be a finite simplicial complex of dimension $n$.
Then there is an snc variety $X$ such that
\begin{enumerate}
\item $\DD(X)\cong \CC$,
\item $\dim X=\dim \CC$,
\item $X$ is a hypersurface in a smooth, projective variety $Y$,
\item $X$ is defined over $\q$ and 
\item the strata of $X$ and also $Y$ are all rational varieties.
\end{enumerate}
\end{prop}

Proof. Let $\{v_i:i\in I\}$ be the vertices of $\CC$.
In $\p^{n+1}$ let $\{H_i:i\in I\}$ be hyperplanes in general
position. 

One can get concrete examples by fixing 
distinct numbers $a_i$ and setting
$$
H_i=\bigl(x_0+a_ix_1+a_i^{2}x_2+\cdots +  a_i^{n+1}x_{n+1}=0\bigr)
\subset \p^{n+1}.
$$
No $n+2$ of the $H_i$ have a point in common since
the coefficients form a  Vandermonde matrix.
If $a_i\in \q$, then $H_i$ is defined over $\q$.
If the $a_i$ are conjugate over $\q$, their union is 
defined over $\q$ but has a nontrivial Galois action on its
geometric irreducible components.

Set $X_0:=\cup_{i\in I} H_i$ and $Y_0:=\p^{n+1}$.

Note that $X_0$ is an snc variety, all of its strata are
linear spaces and for $0\leq r\leq n$ the $(n-r)$-dimensional strata
of $X_0$ are in a natural bijection with 
 $\binom{I}{r}$,  the set of $r$-element subsets of $I$.

Thus every subset $J\subset I$ with at most $n+1$ elements spans
a cell in $\DD(X_0)$ hence 
 there is a natural embedding
$\tau_0:\CC\into \DD(X_0)$ that is  a bijection on vertices.
Thus $X_0$ and $\tau_0$ satisfy the assumptions of
Lemma \ref{snc.bu.induct.lem} with $r=0$
since neither $\CC$ nor $\DD(X_0)$ contain any cells of dimension $>n$.
 
Using Lemma \ref{snc.bu.induct.lem} inductively for $0\leq r< n$,
we obtain  embedded snc varieties
$X_r\subset Y_r$ and 
smooth subvarieties $Z_r\subset X_r$ such that
$$
X_{r+1}=B_{Z_r}X_r\qtq{and}   Y_{r+1}=B_{Z_r}Y_r.
$$
Furthermore, $\tau_0$ restricts to embeddings
$\tau_r:\CC\into \DD(X_r)$
which are   bijections on the vertices
and on cells of dimension $\geq n-r+1$.

Finally we get $X_n$ such that
$\CC\into \DD(X_{n})$  is  a bijection on   cells  in every dimension.

The strata of $X_r$ are obtained from
 linear subspaces in  $Y_0:=\p^{n+1}$ by blowing up lower dimensional
linear subspaces, hence the strata are rational.

Therefore $X:=X_n$ satisfies all the requirements. 
\qed
\medskip

Using the pair $X\subset Y$ obtained in
Proposition \ref{pf.main.main.1}, the following shows
(\ref{main.main.thm}.2).

\begin{lem} \label{log.ample.snc.lem}
Let $Y$ be a smooth, projective variety of dimension $n$,
$D\subset Y$ an snc divisor and
$\{p_j:j\in J\}$ the $0$-dimensional strata. Let $\tau:Y_1\to Y$ be a
double cover ramified along a general, sufficiently  ample
divisor $H'\sim 2H$. For each $p_j$ pick a preimage $q_j\in \tau_1^{-1}(p_j)$
and let $\pi:Y_2\to Y_1$ be the blow-up of the points
$\{q_j:j\in J\}$. Set $D_1:=\tau^*D$ and $D_2:=\pi^{-1}_*D_1$.  Then
\begin{enumerate}
\item $Y_2$ is a smooth, projective variety and 
$D_2\subset Y_2$ is an snc divisor,
\item $K_{Y_2}+D_2$ is ample and
\item the composite $\tau\circ \pi$ induces an isomorphism
$\DD(D_2)\cong \DD(D)$.
\end{enumerate}
\end{lem}

Proof. By the Hurwitz formula,
 $$
K_{Y_1}+D_1\simq \tau^*\bigl(K_Y+D+ H\bigr)
$$
thus if $H$ is sufficiently  ample then 
$K_{Y_1}+D_1 $ is the pull-back of a very ample line bundle
from $Y$. (As for hyperelliptic curves, usually
$K_{Y_1}+D_1 $ is not very ample.)
We may even arrange that
$$
\o_{Y}\bigl(K_{Y}+D+ H\bigr)\bigl(-\tsum_j p_j\bigr)
\eqno{(\ref{log.ample.snc.lem}.4)}
$$
is very ample on $Y\setminus\{p_j:j\in J\}$.

Let $E\subset Y_2$ be the exceptional divisor of $\pi$; it is a union
of $|J|$ disjoint copies of $\p^{n-1}$.
Then
$$
K_{Y_2}\sim \pi^*K_{Y_1}+(n-1)E\qtq{and}
D_2\sim \pi^*D_1-nE.
$$
Therefore
$$
K_{Y_2}+D_2\sim \pi^*\bigl(K_{Y_1}+D_1\bigr)(-E).
$$
This shows that  $K_{Y_2}+D_2 $ is ample on $E$ and
(\ref{log.ample.snc.lem}.4) implies that it is ample on
$Y_2\setminus E\cong Y_1\setminus\{q_j:j\in J\}$.
Thus $K_{Y_2}+D_2$ is ample.

Finally, every positive dimensional stratum  $W_1\subset D_1$ is a
double cover of a stratum  $W\subset D$ ramified along
$W\cap H'\neq \emptyset$.  Thus
$\tau_*:\DD(D_1)\to \DD(D)$ is a bijection on
the cells of dimension $<n$
but every $n$-dimensional cell of  $\DD(D)$ has 2 preimages in 
$\DD(D_1) $. The blow-up of the points $q_j$ removes one
of these cells, giving the isomorphism
$\DD(D_2)\cong \DD(D)$. \qed

\section{Construction of singularities}

Starting with a simplicial complex $\CC$, in
Theorem \ref{main.main.thm} we obtained an embedded
snc variety $D\subset Y$ such that $\DD(D)\cong \CC$. 
If $D\subset Y$ can be contracted to a point  then we get
an isolated singularity $(x\in X)$ one of whose
log resolutions is $(D\subset Y)$. 

In general $D$ is not contractible; we have to perform
some blow-ups first.

\begin{lem} \label{Z.bu.lem}
Let $Y$ be a smooth variety, 
$D\subset Y$ an snc divisor and $H\subset Y$ another
 divisor intersecting $D$ transversally.
 Set $Z:=D\cap H$ and let $\pi:Y_1:=B_ZY\to Y$ be the blow-up with
exceptional divisor $E$. Let $D_1:=\pi^{-1}_*D$ denote the
birational transform of $D$. Then
\begin{enumerate}
\item $D_1\sim \pi^*D-E$ is a Cartier divisor,
\item $D_1\cong D$,
\item $Y_1\setminus D_1$ is smooth,
\item the singularities of  $(Y_1, D_1) $ 
are analytically equivalent to
$$
(x_1\cdots x_n=z=0)\subset (x_1\cdots x_n=yz)\subset \a^{n+m+2},
$$
\item   $(Y_1, D_1) $  is dlt  and
\item $N_{D_1, Y_1}\cong  \o_{Y_1}(D_1)|_{D_1}\cong  \o_Y(D-H)|_D$.
\end{enumerate}
\end{lem}

Proof. The first 5 assertions are local on $Y$, hence we may assume
that there are local coordinates $y_1,\dots, y_n$
such that $D=(y_1\cdots y_r=0)$ and $H=(y_n=0)$ for some $r<n$.
The  blow-up is given by
an equation
$$
\bigl(y_1\cdots y_r\cdot s=y_n\cdot t\bigr)\subset Y\times  \p^1_{st}
\qtq{and} D_1=(y_1\cdots y_r=t=0)=(t=0).
$$
In the affine chart $t\neq 0$, we get the equation
$y_1\cdots y_r\cdot (s/t)=y_n$ which shows that $Y_1\setminus D_1$ is smooth.

The singularities of $(Y_1, D_1)$ are locally of the form
$$
(y_1\cdots y_r=v=0)\subset \bigl(y_1\cdots y_r=y_n\cdot v\bigr)
\subset Y\times  \a^1_{v}.
$$
This shows (4) and (5) is proved in  Lemma \ref{elem.dlt.lem}.

The normal bundle of $D_1$ is $\o_{Y_1}(D_1)|_{D_1}$.
Note that $D_1\sim \pi^*D-E$ and $\pi^{-1}_*H\sim \pi^*H-E$
is disjoint from $D_1$. Thus 
$$
\o_{Y_1}(D_1)\cong \o_{Y_1}\bigl(\pi^*D-\pi^*H\bigr)\cong \pi^*\o_Y(D-H).
$$
Restricting to $D_1\cong D$ gives (6). \qed

\begin{lem}\label{elem.dlt.lem} Let  $Y$ be a normal  variety and
$D\subset Y$  a codimension 1 subvariety. Assume that  
$D$ is  an snc variety and
 at every singular point $x\in \sing Y$
there are local  analytic coordinates such that
 the pair $(Y,D)$ is given as a product of an snc 
 pair with one of the form
$$
\bigl((y_1\cdots y_r=y_{n+1}=0)\subset (y_1\cdots y_r=y_ny_{n+1})\bigr)
\subset \a^{n+1}
$$
for some $0\leq r<n\leq \dim Y$.
Then $(Y,D)$ has a small log resolution
$\pi:(Y^*,D^*)\to (Y,D)$ and $\DD(D^*)=\DD(D)$.
In particular, $(Y,D)$ is dlt.
\end{lem}

Proof. Let $D_0\subset D$ be an irreducible component.
Let $\pi_0:Y'\to Y$ denote the blow-up of $D_0$.
Then $\pi_0$ is an isomorphism at the points where $D_0$ is Cartier.
At other points, in suitable local coordinates we have
$$
(D_0\subset D\subset Y)\cong 
\bigl((y_r=y_{n+1}=0)\subset (y_1\cdots y_r=y_{n+1}=0)\subset 
(y_1\cdots y_r=y_ny_{n+1})\bigr).
$$
The blow-up is covered by 2 coordinate charts.
In one of them we have new coordinate
$y'_r=y_r/y_{n+1}$ and the local equations are
$$
(D'\subset Y')\cong \bigl(( y_{n+1}=0)\subset 
(y_1\cdots y_{r-1}\cdot y'_r=y_n)\bigr).
$$
Thus, in this chart, $D'$ has only 1 irreducible component
and  $(Y', D')$ is an snc pair.

In the other chart the new coordinate is
$y'_{n+1}=y_{n+1}/y_r$ and the local equations are
$$
\bigl( (y_r=0)\subset  (y_r=0)\cup (y_{1}\cdots y_{r-1}=y'_{n+1}=0)
\subset (y_1\cdots y_{r-1}=y_ny'_{n+1})\bigr).
$$
Thus  the $y_r$-coordinate splits off as a direct factor
and the remaining equation is
$$
\bigl(  (y_{1}\cdots y_{r-1}=y'_{n+1}=0)
\subset (y_1\cdots y_{r-1}=y_ny'_{n+1})\bigr).
$$
We also see that the dual complex of $D$ is unchanged.

Performing these blow-ups for every irreducible component of $D$
we end up with $\pi:Y^*\to Y$  such that
$\bigl(Y^*, D^*:=\pi^{-1}_*D\bigr)$ is an snc pair and $\DD(D^*)=\DD(D)$.

Since $\pi$ is small, the discrepancies are unchanged as we go from
$(Y,D)$ to $\bigl(Y^*, D^*\bigr)$ by \cite[2.30]{km-book}.
Since $\pi$ does not contract any stratum of $D^*$, we see that
$(Y,D)$ is dlt.\qed

\begin{rem} \label{res.ordering.rem}
Note that the above resolution process requires an ordering of
the irreducible components of $D$. Thus it does not apply 
if the {\em simple} normal crossing variety $D$ is replaced by a
normal crossing variety. This makes the study of the dual complex
of  resolutions with normal crossing exceptional divisors quite a bit harder.
\end{rem}


\begin{say}[Proof of Theorem \ref{main.sing.thm}]\label{pf.main.sing.thm}
Start with a simplicial complex $\CC$ and use 
Theorem \ref{main.main.thm} to obtain an 
snc pair $D\subset Y$ such that $\DD(D)\cong \CC$. 

Let $H\subset Y$ be a smooth
 divisor intersecting $D$ transversally
such that $H-D$ is ample.
As  in Lemma \ref{Z.bu.lem}, 
  and let $\pi:Y_1\to Y$ be the blow-up of $Z:=D\cap H$
and set  $D_1:=\pi^{-1}_*D$.
By (\ref{Z.bu.lem}.2), $D_1\cong D$ and by (\ref{Z.bu.lem}.6)
its normal bundle is
$N_{D_1, Y_1}\cong  \o_Y(D-H)|_D$. Thus $N_{D_1, Y_1}^{-1}$ is
ample. Therefore, by \cite{artin}, $D_1\subset Y_1$ 
can be contracted to a point, at least analytically or 
 \'etale locally. That is, there is an analytic space
(resp.\ an algebraic variety)  $U_1$ containing $D_1$ 
with an open embedding (resp.\ an \'etale morphism)
$g:(D_1\subset U_1)\to (D_1\subset Y_1)$ 
and a contraction morphism
$$
\begin{array}{ccc}
D_1 & \subset & U_1\\
\downarrow && \hphantom{\pi}\downarrow\pi \\
x & \in & X
\end{array}
\eqno{(\ref{pf.main.sing.thm}.1)}
$$
such that $X$ is normal, $U_1\setminus D_1\cong X\setminus \{x\}$
 and $D_1=\supp\pi^{-1}(x)$.

By (\ref{Z.bu.lem}.4)  $U_1$ is not smooth, but Lemma \ref{elem.dlt.lem}
shows that there is a log resolution 
$$
\tau: (D_2\subset U_2)\to (D_1\subset U_1)
$$
such that $\DD(D_2)=\DD(D_1)$. 
Thus $\pi\circ\tau$ gives a
resolution
$$
\begin{array}{ccl}
D_2 & \subset & U_2\\
\downarrow && \ \downarrow\pi\circ\tau \\
x & \in & X
\end{array}
\eqno{(\ref{pf.main.sing.thm}.2)}
$$
This proves Theorem \ref{main.sing.thm}. \qed
\end{say}

\begin{rem} Although, as we noted, the above construction
does not seem to be compatible with  symmetries of $\CC$,
it is usually possible to realize symmetries  of $\CC$ as Galois
group actions on  $D_1$ (though not on $D_2$ by 
Remark \ref{res.ordering.rem}). This may help in understanding
the algebraic fundamental groups of links.
\end{rem}

\begin{say}[Proof of Theorem \ref{main.rtl.thm}]\label{pf.main.rtl.thm}
The proof is essentially in \cite{k-fg}; see also  \cite[Sec.8]{k-link}.

The implication (\ref{main.rtl.thm}.2) $\Rightarrow$ (\ref{main.rtl.thm}.1)
was established in the cited papers but the key ingredients
are contained already in   
\cite[pp.68--72]{gri-sch} and \cite[2.14]{steenbrink}.

 The converse statements
are not fully proved in \cite[Sec.8]{k-link},
 but the reason  is that instead of our 
Theorem \ref{main.main.thm}, weaker existence results were used.
Using Paragraph \ref{pf.main.sing.thm}, the proof is even simpler.

Let $\CC$ be a $\q$-acyclic simplicial complex.
First, \cite[3.63]{kk-singbook} shows that
with $D$ as in (\ref{main.main.thm}.1) we have
$H^i(D, \o_D)=0$ for $i>0$. Then \cite[3.54]{kk-singbook}
implies that for $H$ sufficiently ample, the singularity
$(x\in X)$ obtained in (\ref{pf.main.sing.thm}.1--2)
is rational.\qed

\end{say}

\begin{say}[Proof of Theorem \ref{main.DMR.thm}]\label{pf.main.DMR.thm}
Start with a simplicial complex $\CC$ and use
Theorem \ref{main.main.thm} to obtain an embedded
snc variety $D\subset Y$ such that $\DD(D)\cong \CC$. 

As before, we blow up a subvariety of the form
$D\cap H$ and then contract the birational transform of
$D$ to obtain the required singularity $(x\in X)$.
However, we have to choose $H$ with some care.

By (\ref{main.main.thm}.2) we may assume that
$K_Y+D$ is ample. Choose $m\gg 1$ such that
$m(K_Y+D)+D$ is very ample and let
$H\sim m(K_Y+D)+D$ be a smooth divisor that intersects $D$ transversally.
Let $\pi:Y_1\to Y$ be the blow-up of $Z:=D\cap H$ with
exceptional divisor $E$.  Then
$$
K_{Y_1}\sim \pi^*K_Y+E,\quad D_1:=\pi^{-1}_*D\sim \pi^*D-E\qtq{and}
\pi^{-1}_*H\sim \pi^*H-E.
$$
Therefore
$$
mK_{Y_1}+(m+1)D_1\sim \pi^*\bigl(mK_Y+(m+1)D\bigr)-E
\sim \pi^{-1}_*H.
\eqno{(\ref{pf.main.DMR.thm}.1)}
$$
Let $(D_1\subset U_1)\to (D_1\subset Y_1)$ be any
analytic or \'etale  neighborhood of $D_1$
whose image is disjoint from $\pi^{-1}_*H $.
Then (\ref{pf.main.DMR.thm}.1) implies that
$$
\o_{U_1}\bigl(mK_{U_1}+(m+1)D_1\bigr)\cong \o_{U_1}.
\eqno{(\ref{pf.main.DMR.thm}.2)}
$$
Let $g:(D_1\subset U_1)\to (x_1\in X_1)$ be the contraction of
$D_1$ to a normal singularity. 
Then the isomorphism (\ref{pf.main.DMR.thm}.2)
pushes forward to an isomorphism
$$
\o_{X_1}(mK_{X_1})\cong g_*\o_{U_1}\bigl(mK_{U_1}+(m+1)D_1\bigr)\cong
g_*\o_{U_1} \cong \o_{X_1}.
\eqno{(\ref{pf.main.DMR.thm}.3)}
$$
Thus  $mK_{X_1}$ is Cartier. 
Furthermore, $(U_1, D_1)$ is dlt by (\ref{Z.bu.lem}.4)
and
$$
\o_{U_1} \bigl(K_{U_1}+D_1\bigr)\bigr|_{D_1}\cong 
\o_{D_1}(K_{D_1})\cong \o_{D}(K_{D})\cong 
\o_{Y} \bigl(K_{Y}+D\bigr)\bigr|_{D}
$$
is ample by assumption, hence 
$g:U_1\to X_1$ is a dlt modification.  Therefore, as noted in
Definition \ref{dres.defn.II}, 
$\DMR(x_1\in X_1)$ is PL-homeomorphic to 
$\DD(D_1)\cong \DD(D)\cong \CC$.
This gives a requisite example where the canonical class is
$\q$-Cartier. 
From this we  get an example where the canonical class is
Cartier as follows.

The isomorphism 
 $\o_{X_1}\bigl(mK_{X_1}\bigr)\cong \o_{X_1}$
 determines a degree $m$ cyclic cover
$$
\tau:(x\in X)\to (x_1\in X_1).
$$ 
(See \cite[2.49--53]{km-book} for the construction and 
basic poperties of cyclic covers.)
Let $\pi:(U_2, D_2)\to (U_1, D_1)$ be the log resolution
constructed in Lemma \ref{elem.dlt.lem} such that
$\DD(D_2)\cong\DD(D_1)$.
Corresponding  to $X\to X_1$ we have a 
cyclic cover $\tau_U:U\to U_2$ that comes from the isomorphism
$$
\o_{U_2}\bigl(K_{U_2}+D_2\bigr)^{\otimes m}\sim \o_{U_2}\bigl(-D_2\bigr).
$$
This shows that $D:=\supp \tau_U^{-1}(D_2)\to D_2$
is an isomorphism. In particular, $U$ and $X$  are irreducible.
 $K_{U}+D $ is not relatively ample any more, but
$\bigl(K_{U}+D\bigr)|_{D}\sim K_{D}\sim K_{D_2}$ is the pull-back of
$K_{D_1}\sim \bigl(K_{U_1}+D_1\bigr)|_{D_1}$, hence
 $K_{U}+D $ is  nef over $X$.

The singularities of  $D\subset U$ are  locally  of the form
$$
\bigl((y_1\cdots y_r=y_n=0)\subset (y_1\cdots y_r=y_n^m)\bigr)
\subset \a^{n+1}
$$
These pairs are not dlt, but
they are quotients of dlt pairs.
Indeed, they can be written as
$$
\bigl((u_1\cdots u_r=0)\subset \a^{n}_{\mathbf u}\bigr)/G_m
$$
where $G_m$ is the subgroup of diagonal matrices
$$
\bigl\{\diag(a_1,\dots, a_r, 1,\dots,1) : 
a_1^m=\cdots =a_r^m=a_1\cdots a_r=1\bigr\}.
$$

Thus, using  \cite[Cor.38]{dkx},
$\CC\cong \DD(D_1)\cong \DD(D_2)\cong \DD(D)$ 
is PL-homeomorphic to  $\DMR(x\in X)$.
\qed

\end{say}

\section{Fundamental groups of rational links}

A group $G$ is called $\q$-superperfect
if $H_1(G,\q)=H_2(G,\q)=0$, see \cite[Def.40]{k-fg}. 
By \cite[Thm.42]{k-fg}, for every finitely presented,
$\q$-superperfect group $G$ there is a 
6-dimensional rational singularity $(x\in X)$
such that $\pi_1\bigl(\operatorname{Lk}(x\in X)\bigr)$,
 the fundamental group of  
 the link of  $(x\in X)$, 
is isomorphic to $G$. Our results yield such examples
also in dimensions 4 and 5. In addition, we get isolated singularities.

\begin{thm}\label{rtl.link.thm}
Let $G$ be a finitely presented,
$\q$-superperfect group. For $n\geq 4$, there is an
$n$-dimensional  rational singularity $(x\in X)$
such that $\pi_1\bigl(\operatorname{Lk}(x\in X)\bigr)\cong G$.
\end{thm}

Proof.  By Lemma \ref{3.dim.top.lem}, there is a
3-dimensional, $\q$-acyclic,  simplicial complex
$\CC$ such that $\pi_1(\CC)\cong G$. 
Thus, for  $n\geq 4$, Theorem \ref{main.rtl.thm} 
constructs an 
$n$-dimensional  rational singularity $(x\in X)$
such that $\pi_1\bigl(\DR(x\in X)\bigr)\cong G$.

By \cite[Thm.35]{k-fg} there is a surjection
$$
\pi_1\bigl(\operatorname{Lk}(x\in X)\bigr)\onto \pi_1\bigl(\DR(x\in X)\bigr)
\eqno{(\ref{rtl.link.thm}.1)}
$$
whose kernel is finite and cyclic.
The discussion at the end of \cite[Sec.6]{k-fg} shows how to modify
the construction to ensure that
(\ref{rtl.link.thm}.1) is an isomorphism.\qed
\medskip

\begin{rem} As discussed in \cite[Sec.7]{k-fg}, if
$(x\in X)$ is a 3-dimensional rational singularity then
$\pi_1\bigl(\DR(x\in X)\bigr) $ has a balanced presentation;
that is, the number of relations equals the number of generators.
For arbitrary rational singularities (\ref{rtl.link.thm}.1)
can have a very large kernel which is very hard to control in general.
 It is not known which groups
can occur as fundamental groups of links of rational singularities
in dimension 3 or higher.
\end{rem}

The following  is well known.

\begin{lem}\label{3.dim.top.lem} Let $G$ be a finitely presented  group.
\begin{enumerate}
\item  There is a 2-dimensional simplicial complex
$\CC_2$ such that $\pi_1(\CC_2)\cong G$.
\item  There is a 3-dimensional simplicial complex
$\CC_3$ such that 
$$
\pi_1(\CC_3)\cong G, \quad  H_i(\CC_3, \z)\cong H_i(G, \z)\ \mbox{for}\
i=1,2
\qtq{and} H_3(\CC_3, \z)\cong 0.
$$
\item  If $G$ is $\q$-superperfect then 
$\CC_3$ is  $\q$-acyclic.
\end{enumerate}.
\end{lem}

Proof. Let $\langle a_1,\dots, a_r | \omega_1, \dots, \omega_s\rangle$
be a presentation of $G$.

Start with $\CC_1$, a bouquet of $r$ circles.
If $\omega_i$ has word length $l_i$, we can kill it in
$\pi_1(\CC_1)$ by attaching an $l_i$-gon to $\CC_1$
appropriately. We can think of the $l_i$-gon
as made up of $l_i$ triangles. This way we have a
2-dimensional $\Delta$-complex $\CC_2$ whose fundamental group is
$G$. Note that $H_2(\CC_2, \z)$ is free.

Let $\beta_1,\dots, \beta_t\in H_2(\CC_2, \z)$ be a basis of the
image of the Hurewich map
$\pi_2(\CC_2)\to H_2(\CC_2, \z)$. We  attach 3-balls
$B_i$ to the $\beta_i$ to obtain a 
3-dimensional $\Delta$-complex $\CC_3$.
By construction $\pi_1(\CC_3)\cong G$ 
(hence also $H_1(\CC_3, \z)\cong H_1(G, \z) $) and
$ H_3(\CC_3, \z)= 0$.

By Hopf's theorem, 
$$
H_2(G, \z)\cong \coker\bigl[\pi_2(\CC_2)\to H_2(\CC_2, \z)\bigr]\cong
H_2(\CC_3, \z)
$$
which completes the proof of  (2). 
If $G$ is $\q$-superperfect then 
$H_i(\CC_3, \q)=0$ for $i=1,2$, hence then 
$\CC_3$ is  $\q$-acyclic.\qed

\begin{ack}
I thank C.~Xu for comments and corrections.
Partial financial support    was provided  by  the NSF under grant number 
DMS-07-58275 and by the Simons Foundation. 
\end{ack}

\def\cprime{$'$} \def\cprime{$'$} \def\cprime{$'$} \def\cprime{$'$}
  \def\cprime{$'$} \def\cprime{$'$} \def\dbar{\leavevmode\hbox to
  0pt{\hskip.2ex \accent"16\hss}d} \def\cprime{$'$} \def\cprime{$'$}
  \def\polhk#1{\setbox0=\hbox{#1}{\ooalign{\hidewidth
  \lower1.5ex\hbox{`}\hidewidth\crcr\unhbox0}}} \def\cprime{$'$}
  \def\cprime{$'$} \def\cprime{$'$} \def\cprime{$'$}
  \def\polhk#1{\setbox0=\hbox{#1}{\ooalign{\hidewidth
  \lower1.5ex\hbox{`}\hidewidth\crcr\unhbox0}}} \def\cdprime{$''$}
  \def\cprime{$'$} \def\cprime{$'$} \def\cprime{$'$} \def\cprime{$'$}
\providecommand{\bysame}{\leavevmode\hbox to3em{\hrulefill}\thinspace}
\providecommand{\MR}{\relax\ifhmode\unskip\space\fi MR }
\providecommand{\MRhref}[2]{%
  \href{http://www.ams.org/mathscinet-getitem?mr=#1}{#2}
}
\providecommand{\href}[2]{#2}

\vskip1cm

\noindent Princeton University, Princeton NJ 08544-1000

{\begin{verbatim}kollar@math.princeton.edu\end{verbatim}}

\end{document}